\documentclass{amsart}

\usepackage{amsfonts}
\usepackage{amsmath}
\usepackage{amsthm}
\usepackage{amssymb}
\usepackage[english]{babel}
\usepackage{graphics}
\usepackage{latexsym}
\usepackage{longtable} \usepackage{mathrsfs}
\usepackage{graphicx}
\usepackage{wrapfig} 
\usepackage[pdftex,colorlinks=true]{hyperref}

\newtheorem{theorem}{Theorem}
\newtheorem{corollary}[theorem]{Corollary}
\newtheorem{lemma}[theorem]{Lemma}
\newtheorem{proposition}[theorem]{Proposition}

\newtheorem{claim}[theorem]{Claim}
\newtheorem{example}[theorem]{Example}

\theoremstyle{definition}
\newtheorem{definition}[theorem]{Definition}

\newcommand{\R}{\mathbb{R}}

\newcommand{\mB}{\mathbb{B}}

\newcommand{\noi}{\noindent}
\newcommand{\ms}{\medskip}

\newcommand{\al}{\alpha}
\newcommand{\be}{\beta}
\newcommand{\ga}{\gamma}
\newcommand{\de}{\delta}
\newcommand{\De}{\Delta}
\newcommand{\e}{\varepsilon}

\newcommand{\Om}{\Omega}

\newcommand{\larrow}{\longrightarrow}
\newcommand{\ot}{\otimes}

\newcommand{\ri}{\rightarrow}

\newcommand{\p}{\partial}
\newcommand{\sub}{\subseteq}

\newcommand{\by}{\times}
\newcommand{\rk}{\textrm{rk}}

\newcommand{\sgn}{\textrm{sgn}}
\newcommand{\ess}{\textrm{ess}}

\newcommand{\dist}{\textrm{dist}}

\newcommand{\Div}{\textrm{Div}}

\newcommand{\spn}{\textrm{span}}

\newcommand{\bt}{\begin{theorem}}\newcommand{\et}{\end{theorem}}
\newcommand{\bd}{\begin{definition}}\newcommand{\ed}{\end{definition}}
\newcommand{\bl}{\begin{lemma}}\newcommand{\el}{\end{lemma}}
\newcommand{\beq}{\begin{equation}}\newcommand{\eeq}{\end{equation}}
\newcommand{\bc}{\begin{claim}}\newcommand{\ec}{\end{claim}}
\newcommand{\bex}{\begin{example}}\newcommand{\eex}{\end{example}}
\newcommand{\bcor}{\begin{corollary}}\newcommand{\ecor}{\end{corollary}}
\newcommand{\bp}{\begin{proof}}\newcommand{\ep}{\end{proof}}

\newcommand{\BPL}{\medskip \noindent \textbf{Proof of Lemma} }

\newcommand{\BPP}{\medskip \noindent \textbf{Proof of Proposition} }

\numberwithin{equation}{section}
\numberwithin{theorem}{section}

\begin{document}

\title{The Subelliptic $\infty$-Laplace System on Carnot-Carath\'eodory Spaces}

\author{\textsl{Nicholas Katzourakis}}
\address{Department of Mathematics and Statistics, University of Reading, Whiteknights, PO Box 220, Reading RG6 6AX, Berkshire, UK and Basque Center for Applied Mathematics, Alameda de Mazarredo 14, E48009, Bilbao, Spain.}
\email{n.katzourakis@reading.ac.uk}


\date{}



\begin{abstract} Given a  Carnot-Carath\'eodory space $\Om \sub \R^n$ with associated vector fields $X=\{X_1,...,X_m\}$, we derive the subelliptic $\infty$-Laplace system for mappings $u : \Om \larrow \R^N$, which reads
\[  \label{1}
\De^X_\infty u  \, :=\, \Big(Xu \ot Xu  + \|Xu\|^2 [Xu]^\bot \! \ot I \Big) : XX u\, = \, 0  \tag{1}
\]
in the limit of the subelliptic $p$-Laplacian as $p\ri \infty$. Here $Xu$ is the horizontal gradient and $[Xu]^\bot$ is the projection on its nullspace. Next, we identify the Variational Principle characterizing \eqref{1}, which is the ``Euler-Lagrange PDE'' of the supremal functional 
\[ \label{2}
E_\infty(u,\Om)\ :=  \ \|Xu\|_{L^\infty(\Om)}   \tag{2}
\]
for an appropriately defined notion of \emph{Horizontally $\infty$-Minimal Mappings}. We also establish a maximum principle for $\|Xu\|$ for solutions to \eqref{1}. These results extend previous work of the author \cite{K1, K2} on vector-valued Calculus of Variations in $L^\infty$ from the Euclidean to the subelliptic setting. 
\end{abstract}

\maketitle

\section{Introduction} \label{section1}

Let $X:=\{X_1,...,X_m\}$ with $X_i : \Om\sub \R^n \larrow \R^n$ be a frame of $C^1(\Om)^n$ vector fields defined on the connected domain $\Om \sub \R^n $ with $1\leq m \leq n$. The linear span $H(x) := \spn[\{X_1(x),...,X_m(x)\}]$ is called the Horizontal subspace of $\R^n$ at $x \in \Om$. We equip $\Om$ with a Riemannian metric $g$ that makes $X$ an orthonormal family, that is $g(X_i,X_j)=\de_{ij}$, $1\leq i,j \leq m$. Given a map $u=(u_1,...,u_N)^\top =u_\al e_\al : \Om \sub \R^n \larrow \R^N$, we denote the directional differentiation of $u$ along $X_i$ by $X_i u$ and we define its Horizontal gradient as
\[
Xu := (X_iu_\al) e_\al \ot X_i\ :  \ \Om \sub \R^n \larrow \R^{N \by n}.
\]
Note that $Xu(x) \in \R^N \ot H(x) \sub \R^{N \by n}$. Throughout this paper, the summation convention is employed in repeated indices in a product. Greek indices $\al, \be, \ga,... $ will run from $1$ to $N$, lowercase Latin $i,j,k,...$ from $1$ to $m$ and uppercase Latin $A,B,C,...$ from $1$ to $n$.  We equip $ \R^{N \by n}$ with the induced natural metric for which $g(Xu,Xv)=X_iu_\al X_i v_\al$ and set $\|Xu\|^2:=g(Xu,Xu)$. 

In this paper we are interested in vector-valued Calculus of Variations in the space $L^\infty$ for the model supremal functional
\beq \label{eq1}
E_\infty(u,\Om)\ := \ \big\|Xu\big\|_{L^\infty(\Om)}
\eeq
which we interpret as $\ess \sup_{\Om}\|Xu\|$, and its associated ``Euler-Lagrange PDE system" which we call \emph{subelliptic $\infty$-Laplacian}:
\beq   \label{eq2}
X_i u_\al\, X_ju_\be \,X_iX_j u_\be \ +\ \|Xu\|^2 [Xu]_{\al \be}^\bot X_iX_i u_\be\ = \ 0.
\eeq
Here $[Xu]^\bot(x)$ is the projection on the nullspace of the linear operator $Xu(x)^\top : \R^N \larrow H(x)$. The operator $X_iX_i$ is the well known ``H\"ormander's sum of squares" in subelliptic theory. In compact vector notation, we write \eqref{eq2} as
\beq  \label{eq3}
\De^X_\infty u \, :=\, \Big(Xu \ot Xu  + \|Xu\|^2 [Xu]^\bot \! \ot I \Big) : XX u\, = \, 0. 
\eeq
The Euclidean case of \eqref{eq3} is obtained by taking $m=n$ and as $\{X_1,...,X_n\}$ the partial derivatives $\{D_1,...,D_n\}$. Then we deduce
\beq  \label{eq4}
\De_\infty u \, =\, \Big(Du \ot Du  + |Du|^2 [Du]^\bot \! \ot I \Big) : D^2 u\, = \, 0
\eeq
where the respective functional is the $L^\infty$ norm of the Euclidean norm on $\R^{N \by n}$ of the gradient $Du$, i.e.\ $|Du|=(D_iu_\al D_iu_\al)^{\frac{1}{2}}$. System \eqref{eq4} has first been derived by the author in \cite{K1} and has been subsequently studied alongside its associated functional in \cite{K2,K3}. System \eqref{eq4} is a quasilinear degenerate elliptic system in non-divergence form (with discontinuous coefficients) which arises in the limit of the $p$-Laplace system $\De_p u = \Div \big(|Du|^{p-2}Du\big)=0$ as $p\ri \infty$. The special case of the scalar $\infty$-Laplace PDE for $N=1$ reads
$\De_\infty u = D_iu\, D_ju \,D^2_{ij}u=0$ and has a long history. In this case the coefficient $|Du|^2[Du]^\bot$ of \eqref{eq4} vanishes identically. The scalar $\De_\infty$ was derived in the limit of the $p$-Laplacian as $p\ri \infty$ in the '60s by Aronsson and was first studied in \cite{A3, A4} (see also \cite{A1,A2}). It has been extensively studied ever since, but most of the associated problems have been solved in the last 20 years in the context of Viscosity Solutions (see for example Crandall \cite{C} and references therein). 

A basic difficulty associated to \eqref{eq4} which is a genuinely vectorial phenomenon and does not appear in the scalar case is that $|Du|^2[Du]^\bot$ may be discontinuous even for $C^\infty$ solutions. Such an example on $\R^2$ is given by $u(x,y) = e^{ix}-e^{iy}$, which is $\infty$-Harmonic near the origin but the projection $[Du]^\bot$ is discontinuous on $\{x=y\}$, since the rank of $Du$ jumps from 2 to 1  on the diagonal. In general, $\infty$-Harmonic maps present a phase separation,  which is quite well understood in two dimensions (\cite{K1,K3}). Much more intricate examples of smooth $\infty$-Harmonic maps in two dimensions whose interfaces are not straight lines but instead have triple junctions and corners are constructed in the very recent paper \cite{K5}. 

The motivation to study $L^\infty$ variational problems stems from their frequent appearance in applications (see e.g.\ \cite{B} for the scalar case) because minimizing maximum values (e.g.\ maximum tensions before fraction) furnishes more realistic models when compared to minimization of averages which corresponds to integral functionals. They also are analytically extremely interesting and in particular the related equations are in nondivergence form and with discontinuous coefficients. Moreover, certain geometric problems are inherently connected to $L^\infty$. In the vector case $N\geq 2$ our motivation comes from the problem of optimization of quasiconformal deformations of Geometric Analysis (see \cite{CR} and \cite{K4}). In the scalar case $N=1$, the main motivation come from the optimization of Lipschitz Extensions (\cite{A3,C}).

From the variational viewpoint, a central difficulty arising in the study of \eqref{eq1} is that it is \emph{nonlocal}, in the sense that with respect to the $\Om$ argument it is not a measure. This implies that minimizers over a domain with fixed boundary values are not local minimizers over subdomains and the direct method of Calculus of Variations when applied to \eqref{eq1} does not produce PDE solutions of \eqref{eq2}. In \cite{K2} we identified the appropriate variational notion governing $\infty$-Harmonic maps $u:\Om \sub \R^n \larrow \R^N$. We introduced the concept of \emph{$\infty$-Minimal Maps}, which is a weak version of minimizer with respect to two sets of local variations (reflecting the non-divergence form of \eqref{eq4}): essentially scalar local variations with fixed boundary values (``Rank-One Absolute Minimality")  and normal free variations (``$\infty$-Minimal Area of the submanifold $u(\Om) \sub \R^N$"). 

Herein, following \cite{K1}, we derive \eqref{eq2} in the limit of the subelliptic Euler-Lagrange equation of the Horizontal $p$-Dirichlet functional, that is the $L^p$-norm of the Horizontal gradient $\|Xu\|_{L^p(\Om)}$. Observe that at least in a formal level we have
\[
\De^X_p u \, \larrow \, \De^X_\infty u \ \ \text{ and }\ \  \|Xu\|_{L^p(\Om)}\, \larrow \, \|Xu\|_{L^\infty(\Om)},
\]
both as $p\ri \infty$, but it is not a priori clear that the following rectagle ``commutes" 
\begin{align} \label{1.3}
&\|Xu\|_{L^p(\Om)}   \ \ \ \  \larrow\ \ \ \ \De^X_p u =0  \nonumber\\
& \ \ \ \downarrow\ p \ri \infty\ \ \ \ \ \ \ \ \ \ \  \ \ \ \ \ \ \downarrow\ p \ri \infty  \nonumber\\
&\|Xu\|_{L^\infty(\Om)}   \ \ \ \ \dashrightarrow \ \ \ \ \De^X_\infty u=0   \nonumber
\end{align}
in the sense that \eqref{eq1} and \eqref{eq3} are directly related. Next, inspired by \cite{K2}, we introduce a subelliptic variant of $\infty$-Minimal Maps which we call \emph{Horizontally $\infty$-Minimal maps}, and establish equivalence between these local minimizers of \eqref{eq1} and solutions of \eqref{eq2}. Interestingly,  the variational problem for \eqref{eq1} is sufficient for \eqref{eq2} for an arbitrary $m$-frame of vector fields $\{X_1,...,X_m\}$ on $\Om$ without extra assumptions, but sufficiency holds only when $\Om$ is a \emph{Carnot-Carath\'eodory metric space} with respect to $\{X_1,...,X_m\}$ (see definition below). Finally, we also introduce an Horizontal gradient flow associated to \eqref{eq2} and by using that tool we establish Maximum-Minimum Principles for $\|Xu\|$. The referee of this paper pointed out that our technique in the proof can be characterized as a technique of ``propagation of maxima and minima" along integral curves of vector fields and is analogous to the proof of Bony in \cite{Bo} in the case of smooth H\"ormander vector fields.

In order for our analysis to be made rigorous and precise and focus on the new structures that emerge, we restrict ourselves to the class of \emph{Horizontally $C^1$ maps of full rank}. This class consists of maps for which $Xu$ is continuous and its rank satisfies $\rk(Xu)\equiv \min\{m,N\}$ on $\Om$. Since $Xu(x) \in \R^N \ot H(x)$, generally $\rk(Xu)\leq \min\{m,N\}$. If $Xu=Du$ and $m=n$, this class consists of immersions and submersions. We also impose the extra simplifying assumption that $Xu$ is differentiable in the Euclidean sense, which is not necessary but allows to avoid technical difficulties and regularizations. The restriction on the rank of $Xu$ owes to that the projection coefficient $[Xu]^\bot$ becomes discontinuous when the rank of $Xu$ varies (see \cite{K1,K2,K3} for related analysis of the elliptic version \eqref{eq4}). 

We conclude with some known results related to this paper which inspired and motivated our analysis. In \cite{BC}, Bieske and Capogna studied an extension of the scalar $\infty$-Laplacian called Aronsson's PDE in Carnot groups. In \cite{W}, Wang studied Aronsson's PDE in Carnot-Carath\'eodory spaces arising from vector fields satisfying H\"ormander's condition and subsequently the results have been sharpened by Wang and Yu in \cite{WY}. In the vector case, system \eqref{eq4} has first been derived and studied in \cite{K1}, its variational structure was studied in \cite{K2} and some deeper analytic properties of smooth solutions where studied in \cite{K3}.

\section{Formal Derivation of the Subelliptic $\infty$-Laplace System.} \label{section2}

For an horizontally smooth map $u: \Om \sub \R^n \larrow \R^N$, consider the $p$-Dirichlet functional
\[
E_p(u,\Om)\ := \ \int_\Om\|Xu\|^p
\]
where $Xu$ is the Horizontal gradient, $\|Xu\|^2=g(Xu,Xu)$ and recall that for this Riemannian metric for which $X=\{X_1,...,X_m\}$ on $\Om$ is orthonormal framily we have $\|Xu\|^2=X_i u_\al X_i u_\al$ when expanding $Xu(x)$ as $(X_iu_\al(x)) e_\al\ot X_i(x)$. Recall also that $X_iu(x)=\frac{d}{dt}\big|_{t=0}u(a(t))$ for a curve $a : (-\e,\e)\sub \R \larrow \Om$ with $a(0)=x$ and $a'(0)=X_i(x)$ which can be always chosen to be the affine $a(t)=x+tX_i(x)$. The Euler-Lagrange subelliptic $p$-Laplace PDE system associated to $E_p$ is
\beq \label{eq2.1}
X^*_i\Big(\|Xu\|^{p-2}X_iu_\al\Big)\, =\, 0 
\eeq
where $X^*_i$ is the adjoint differential operator of $X_i$, defined as $X^*_i v:=D_A(X_{iA}v) $
when we expand $X_i(x) =X_{iA}(x)D_A$ in the standard basis of partial derivatives. As usually, we may identify $D_A=\p / \p x_A$ with $e_A=(0,...,1,...,0)^\top$.  For completeness, let us quickly derive \eqref{eq2.1}. Indeed, for $\phi \in C^\infty_0(\Om)^N$ we have
\begin{align}
\frac{d}{dt}\Big|_{t=0}E_p(u+t\phi, \Om)\ &=\ \frac{d}{dt}\Big|_{t=0}\int_\Om \big(X_i (u_\al +t\phi_\al) X_i (u_\al +t\phi_\al)\big)^{p/2} \nonumber\\
&=\ p \int_\Om\|Xu\|^{p-2} X_i u_\al X_i \phi_\al\\
&=-p \int_\Om D_A\Big(X_{iA} \|Xu\|^{p-2} X_i u_\al \Big) \phi_\al\nonumber\\
&=-p \int_\Om X^*_i\Big(\|Xu\|^{p-2}X_iu_\al\Big)\phi_\al. \nonumber
\end{align}
By distributing derivatives in \eqref{eq2.1}, we have
\begin{align}
0\ &=\ D_A\Big(X_{iA} \big(X_k u_\ga X_k u_\ga \big)^{\frac{p-2}{2}}X_iu_\al  \Big)\nonumber\\
    &=\ (p-2) \big(X_k u_\ga X_k u_\ga \big)^{\frac{p-4}{2}}    \Big(X_iu_\al \, X_j u_\be \, X_{iA} D_A X_j u_\be \Big)   \\
&\ \ \ \ +\ \big(X_k u_\ga X_k u_\ga \big)^{\frac{p-2}{2}} \Big[X_i u_\al D_A X_{iA} + X_{iA}\, D_AX_i u_\al\Big].\nonumber
\end{align}
By normalizing and using that $X_i=X_{iA}D_A$, we have
\beq \label{eq2.4}
X_iu_\al \, X_j u_\be \, X_iX_j u_\be  \ +\ \frac{\|Xu\|^2}{p-2} \Big[X_i u_\al D_A X_{iA} + X_iX_i u_\al\Big]\ = \ 0.
\eeq
We define the orthogonal projections of $\R^N$ on the range of the linear map $Xu(x) : H(x)\sub \R^n \larrow \R^N$ and on the nullspace of its transpose $Xu(x)^\top$:
\begin{align}
[Xu]^\top(x)\ &:=\ Proj_{R(Xu(x))},  \label{eq2.5}\\
[Xu]^\bot(x)\ &:=\ Proj_{N(Xu(x)^\top)}. \label{eq2.6}
\end{align}
Note that we have the splitting $I=[Xu]^\top + [Xu]^\bot$ for the identity map of $\R^N$, everywhere on $\Om$. By expanding the term in bracket in \eqref{eq2.4} and observing that $ [Xu]^\bot_{\al \be} X_i u_\be=0$, we have
\begin{align} \label{eq2.7}
X_iu_\al \, X_j u_\be \, X_iX_j u_\be  \ +&\ \frac{\|Xu\|^2}{p-2} [Xu]^\top_{\al \be}\Big[X_i u_\be D_A X_{iA} + X_iX_i u_\be\Big] \nonumber\\
&= -\frac{\|Xu\|^2}{p-2} [Xu]^\bot_{\al \be}\Big[X_i u_\be D_A X_{iA} +  X_iX_i u_\be\Big]\\
&= - \frac{\|Xu\|^2}{p-2} [Xu]^\bot_{\al \be}  X_iX_i u_\be \nonumber.
\end{align}
The crucial observation now is that since
\[
X_iu_\al \, X_j u_\be \, X_iX_j u_\be  \ =\ [Xu]^\top_{\al \ga}X_iu_\ga \, X_i \Big(\frac{1}{2}X_j u_\be X_j u_\be\Big)
\]
the two ends of \eqref{eq2.7} are normal to each other and yet equal. Hence, they both vanish. We choose to multiply the last term $ \frac{\|Xu\|^2}{p-2} [Xu]^\bot_{\al \be}  X_iX_i u_\be$ by $p-2$ and we obtain
\begin{align}
X_iu_\al \, X_j u_\be \, X_iX_j u_\be  \ +&\ \frac{\|Xu\|^2}{p-2} [Xu]^\top_{\al \be}\Big[X_i u_\be D_A X_{iA} + X_iX_i u_\be\Big] \\ 
&=- \|Xu\|^2 [Xu]^\bot_{\al \be}  X_iX_i u_\be.\nonumber
\end{align}
Rearranging, we get
\begin{align}
X_iu_\al \, X_j u_\be \, X_iX_j u_\be \ +&\ \|Xu\|^2 [Xu]^\bot_{\al \be}  X_iX_i u_\be \\ 
&= - \frac{\|Xu\|^2}{p-2} [Xu]^\top_{\al \be}\Big[X_i u_\be D_A X_{iA} + X_iX_i u_\be\Big]\nonumber
\end{align}
and as $p\ri \infty$ we obtain the subelliptic $\infty$-Laplacian \eqref{eq2}.  We note that when $\|Xu\|>0$ the system $\|Xu\|^2 [Xu]^\bot  XX u= 0$ is equivalent to $[Xu]^\bot  XX u= 0$, but we keep the positive function $\|Xu\|^2$ because for ``singular solutions'' these systems generally are not equivalent (cf.\ \cite{K1,K2,K3}).

\section{Variational Characterization of the subelliptic $\infty$-Laplace system.} \label{section3}

We begin with some basics. An absolutely continuous curve $r :[0,T] \larrow \Om$ is called \emph{admissible} when it solves the differential inclusion 
\[
r'(t) \in H(r(t))\ , \ \ a.e.\ t \in [0,T]
\]
and $\|r'\|^2=g(r',r')\leq 1$ a.e.\ on $[0,T]$. This means that there exist measurable coefficients $a_i : [0,T]\larrow \R$ such that
\[
r'(t)\, = \, a_i(t)X_i(t)\  \text{ and }  \ a_i(t)a_i(t) \leq 1\ , \ \ a.e.\ t \in [0,T]
\]  
(the summation convention is employed). For any $x,y \in \Om$, we define the  \emph{Carnot-Carath\'eodory} function
\[
d_X(x,y)\, := \, \inf \Big\{ T \Big| \exists \text{ admissible curve } r :[0,T] \larrow \Om \text{ with } r(0)=x,\ r(T)=y \Big\}.
\]
The domain $\Om$ is called \emph{Carnot-Carath\'eodory space} when $d_X$ is a metric distance on $\Om$, that is when any two points in $\Om$ can be connected by some admissible curve lying into $\Om$. A particular 
class of such spaces are those generated by $m$ smooth vector fields $X_1,...,X_m$ on $\Om$ which sayisfy H\"ormander's condition, that is at each point of $\Om$, the vector fields together with their commutators
\[
X_i,...,\ [X_i,X_j],...,\ [[X_i,X_j],X_k],...
\]
up to a finite order span $\R^n$ everywhere. In this case, the results of \cite{NSW} imply that this integrability condition always guarrantees the existence of admissible curves with horizontal velocities. This fact was first proved by Carath\'eodory in \cite{Ca}. For further material on the subelliptic theory we refer to the papers \cite{H,RS}.

In this section we will consider \eqref{eq1} and \eqref{eq3} in the subclass of horizontally $C^1$ maps (which means $u,X_i u \in C^0(\Om)^N$), for which the horizontal gradient $Xu$ is differentiable in the Euclidean sense:
\[
\Gamma^2(\Om)^N\ :=\ \Big\{ u : \ \Om \sub \R^n \larrow \R^N \, \big| \, u,\, X_iu \in C^0(\Om)^N,\ D(X_iu) \in C^0(\Om)^{N \by n}\Big\}. 
 \]
The latter is a simplifying assumption, replacing the more natural $X_iX_j u \in C^0(\Om)^N$ and allows to bypass regularization schemes. We will also need to impose further restriction on the rank of $Xu$, since the projection coefficient $[Xu]^\bot$ of \eqref{eq3} is discontinuous when the rank of $Xu$ varies. We will consider the class of Horizontal immersions
\[
\Gamma^2_{im}(\Om)^N\ :=\ \Big\{ u  \in \Gamma^2(\Om)^N \, \big| \, \rk(Xu)\equiv m \leq N\Big\},
 \]
the class of Horizontal submersions
\[
\Gamma^2_{sb}(\Om)^N\ :=\ \Big\{ u  \in \Gamma^2(\Om)^N  \, \big| \,  \rk(Xu)\equiv N \leq m\Big\}
 \]
and their union, which is the class of maps with Horizontally full rank:
\[
\Gamma^2_{rk}(\Om)^N\ :=\ \Gamma^2_{im}(\Om)^N \bigcup \Gamma^2_{sb}(\Om)^N\ =\ \Big\{ u  \in \Gamma^2(\Om)^N  \, \big| \,  \rk(Xu)\equiv \min\{m,N\} \Big\}.
 \]
Since $Xu(x) \in \R^N \ot H(x)$, note that generally $\rk(Xu)\leq \min \{m, N\}$. The $\Gamma^1$-variants of these classes are defined by dropping the assumption $D(X_i u) \in C^0(\Om)^{N \by n}$.

\ms

Following \cite{K2}, we introduce a minimality notion for functional \eqref{eq1}. We need the preliminary notion of local \emph{vertical vector fields} to maps $u : \Om \sub \R^n \larrow \R^N$ in $\Gamma^1_{im}(\Om)^N$, relative to the splitting \eqref{eq2.5}, \eqref{eq2.6}: a vertical vector field over $D\sub \Om$ is a smooth map  $\nu : D \sub \Om  \larrow \R^N$ such that $\nu_\be [Xu]_{\be \al}^\top=0$ on $D$. This means $\nu(x) \in [Xu(x)]^\bot$ for all $x\in D$ and consequently 
\beq \label{eq3.1}
\nu_\al\, X_iu_\al = 0.
\eeq
We denote the set of  vertical vector fields over $D\sub \Om$ by $\Gamma([Xu]^\bot,D)$. Note that if the rank of $Xu$ is not constant on $D$, then $\Gamma([Xu]^\bot,D)$ may be empty. 

\begin{definition}\label{def1} Let $u : \Om \sub \R^n \larrow \R^N$ be a map in $\Gamma^1(\Om)^N$.

\ms \noi (i) The map $u$ is called \emph{Horizontal Rank-One Absolute Minimal} on $\Om$ when for all compactly contained subdomains $D$ of $\Om$, all functions $g$ on $D$ vanishing on $\p D$ and all unit  directions $\xi \in \R^N$, $u$ is a minimizer on $D$ with respect to essentially scalar variations $u+g\xi$:
\beq \label{2.2}
\left.
\begin{array}{l}
D \Subset \Om, \\
g\in \Gamma^1_0(D), \\
\xi \in \R^N,\, |\xi|=1
\end{array}
\right\} \ \ \Longrightarrow \ \
E_\infty(u,D)\ \leq \ E_\infty(u+ g\xi , D).
\eeq
\[
\underset{\text{Figure 1.}}{\includegraphics[scale=0.18]{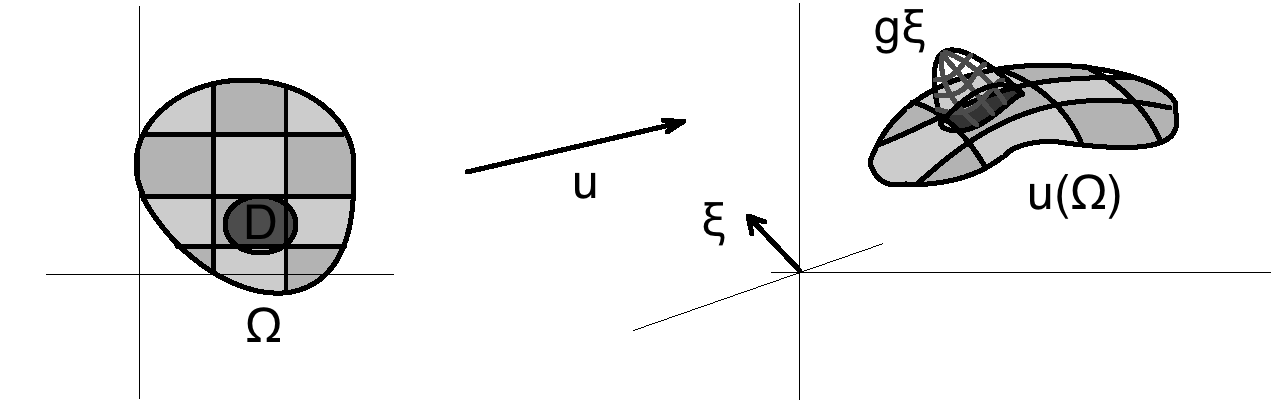}}
\]
In \eqref{2.2} $|\cdot|$ denotes the Euclidean norm of $\R^N$.

\ms

(ii) Suppose further that $\Gamma^1_{rk}(\Om)^N$. We say that \emph{$u(\Om)$ has Horizontally $\infty$-Minimal Area} when for all compactly contained subdomains $D$, all functions $h$ on $\bar{D}$ (not only vanishing on $\p D$) and all vertical vector fields $\nu$, $u$ is a minimizer on $D$ with respect to vertical free variations $u+h\nu$:
\beq \label{2.3}
\left.
\begin{array}{l}
D \Subset \Om, \\
h\in \Gamma^1(\bar{D}), \\
\nu \in \Gamma([Xu]^\bot,D)
\end{array}
\right\} \ \ \Longrightarrow \ \
E_\infty(u,D)\ \leq \ E_\infty(u+ h \nu, D).
\eeq
\[
\underset{\text{Figure 2.}}{\includegraphics[scale=0.18]{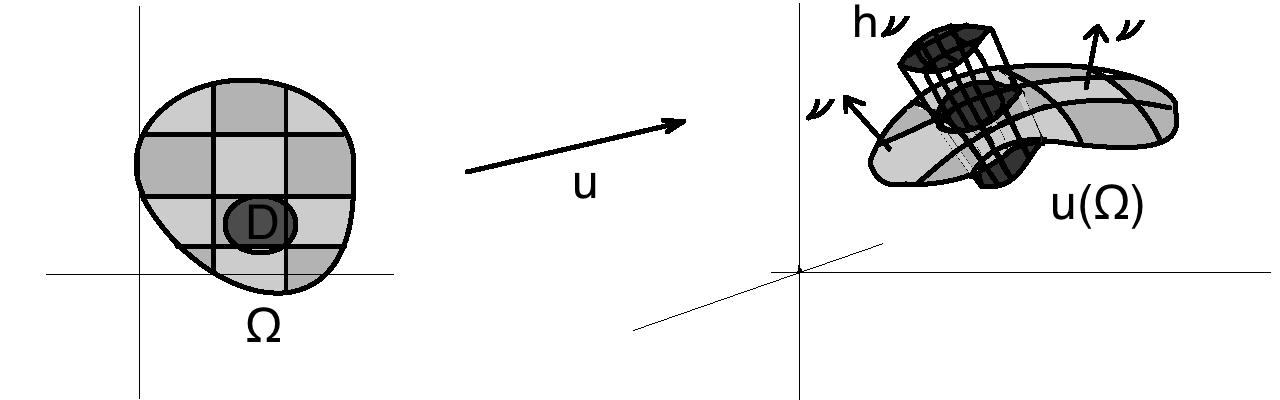}}
\]
(iii) We call $u$ an \emph{Horizontally $\infty$-Minimal Map with respect to the functional \eqref{eq1}} when \eqref{2.2} and \eqref{2.3} hold.
\end{definition}

The main result of this paper is the next

\bt[\textbf{Variational Structure of the subelliptic $\infty$-Laplace system on  Carnot-Carath\'eodory spaces}]

\label{th1}

Let $u : \Om \sub \R^n \larrow \R^N$ be a map with horizontally full rank in $\Gamma^2_{rk}(\Om)^N$ with respect to a frame of vector fields $X=\{X_1,...,X_m\}$ with $X_i \in C^1(\Om)^n$ and $m\leq n$. Suppose also that $\Om$ is a connected open set which is equipped with a Riemannian metric $g$ that makes $X$ orthonormal frame.

\noi Then:

If $u$ is an Horizontally $\infty$-Minimal Map with respect to the functional $E_\infty(u,\Om) = \|Xu\|_{L^\infty(\Om)}  $, then $u$ solves the subelliptic $\infty$-Laplacian 
\[
\De^X_\infty u \ =\, \Big(Xu \ot Xu  + \|Xu\|^2 [Xu]^\bot \! \ot I \Big) : XX u\ = \ 0. 
\]

Conversely, if $\Om$ is a Carnot-Carath\'eodory space with respect to $X=\{X_1,...,X_m\}$ and $u \in \Gamma^2_{im}(\Om)^N$ is an horizontally smooth immersion, then solutions of the subelliptic $\infty$-Laplacian are Horizontally $\infty$-Minimal maps with respect to $E_\infty$.

\et

The proof of Theorem \ref{th1} is split in four lemmas, in each or which the assumptions of Theorem \ref{th1} are subsumed. A basic observation is that the subelliptic $\infty$-Laplace vectorial operator $\De_\infty^X$ splits to two terms, each one normal to the other: since $X_i u_\al X_i(\frac{1}{2}\|Xu\|^2)$ is normal to $[Xu]^\bot_{\al \be}X_iX_iu_\be$, we have that $\De_\infty^Xu =0 $ if and only if 
\begin{align}
X_i u_\al \, X_ju_\be \, X_iX_ju_\be \ &=\  0,\\
\|Xu\|^2[Xu]^\bot_{\al \be}X_iX_iu_\be\ &=\ 0 .
\end{align}
We first have

\bl \label{l1} If $u$ is an Horizontal Rank-One Absolute Minimal on $\Om$, then $u$ solves $X_i u_\al \, X_j u_\be \, X_i X_j u_\be =0$ on $\Om$.
\el

In Lemma \ref{l1} no rank assumption for $Xu$ is needed.

\ms

\BPL \ref{l1}. Fix $x\in \Om$, $0<\e < \dist(x,\p \Om)$, $0<\de<1$ and $\xi \in \R^N$, $|\xi|=1$. Choose $D:=\mB_\e(x)$ the Euclidean $\e$-ball at $x$ and define
\[
 g(z)\ :=\ \frac{\de}{2}\big(\e^2-|z-x|^2\big) \ \in \ \Gamma^1_0(D)\cap\Gamma^2(D).
\]
Set  also $w:=u+g\xi$. Then, by (Euclidean) Taylor expansions of $\|Xu\|^2$ and $\|Xw\|^2$ at $x$ we have
\begin{align} \label{2.7}
\|Xu(z)\|^2\, =\ \|Xu (x)\|^2 \ +\ D_A\big(\|Xu\|^2)(x)(z-x)_A\ + \ o(|z-x|),
\end{align}
as $z\ri x$, and also
\begin{align} \label{2.7a}
\|Xw(z)\|^2\, =\ \|Xw(x)\|^2 \ +\ D_A\big(X_iw_\al X_iw_\al )(x)(z-x)_A\ + \ o(|z-x|),
\end{align}
as $z\ri x$. Using that $X_i = X_{iA}D_A$, we calculate
\beq
X_iw_\al(z) \ =\  X_{iB}(z)D_B\big( u_\al + \xi_\al g \big)(z) \ =\  X_iu_\al(z)- \de \xi_\al X_{iB}(z)(z-x)_B ,
\eeq
and also
\beq
D_A X_i  w_\al(z) \ =\ D_AX_i u_\al(z)-\de \xi_\al \Big[X_{iA}(z)+ D_AX_{iB}(z)(z-x)_B \Big].
\eeq
Consequently, 
\begin{align}
X_iw_\al(x)\ &=\ X_iu_\al(x),\\
D_A X_i  w_\al(x)\ &=\ D_A X_i  u_\al(x) - \de \xi_al X_{iA}(x).
\end{align}
Combining these, we obtain
\begin{align} 
\|Xw(z)\|^2\, &=\ \|Xw(x)\|^2\ +\ 2  X_iw_\al (x)(D_A X_i w_\al)(x) (z-x)_A\ + \ o(|z-x|),\nonumber
\\
             & = \ \|Xu(x)\|^2\ +\ 2 X_iu_\al (x) \Big[  D_A X_i u_\al - \de \xi_\al X_{iA}\Big] (x)(z-x)_A \nonumber\\
     &\ \ \ \  + \ o(|z-x|)\\
 & = \ \ \|Xu(x)\|^2\ +\ \Big[D_A \big( X_iu_\al X_i u_\al \big) - 2\de \xi_\al X_iu_\al X_{iA}\Big](x)  (z-x)_A \nonumber\\
     &\ \ \ \  + \ o(|z-x|),\nonumber
\end{align}
as $z\ri x$. 
Hence,
\begin{align} \label{2.8}
\|Xw(z)\|^2\,  & = \ \|Xu(x)\|^2\ +\ \Big[D_A \big(X_iu_\al X_iu_\al  \big)(x) \nonumber\\
&\ \ \ \  - 2\de \xi_\al \big(X_iu_\al X_{iA}\big)(x) \Big] (z-x)_A  \  + \ o(|z-x|)
\end{align}
as $z\ri x$. 

By \eqref{2.7} we have the estimate
\begin{align}   \label{2.9a}
\left(E_\infty\big(u,\mB_\e(x)\big)\right)^2\ & =\ \sup_{|z-x|<\e} \|Xu(z)\|^2  \nonumber\\
& \geq  \ \|Xu(x)\|^2 \ +\ \max_{\{|z-x|\leq \e\}}\Big\{D_A\big(X_iu_\al X_iu_\al \big)(x) (z-x)_A\Big\}\  \\
&\ \ \ \ \ + \ o(\e), \nonumber
\end{align}
as $\e \ri 0$. By observing that the element $z-x$ which maximizes in \eqref{2.9a} is $\e$ times the sign of $D\big(X_iu_\al X_iu_\al \big)(x)$, we get 
\begin{align} \label{2.9}
\left(E_\infty\big(u,\mB_\e(x)\big)\right)^2\ \geq \ \|Xu(x)\|^2\ +\ \e\big|D\big(X_iu_\al X_iu_\al \big)(x)\big|\ + \ o(\e),
\end{align}
as $\e \ri 0$, where $| \cdot |$ denotes the Euclidean norm on $\R^n$.  Also by \eqref{2.8} we have
\begin{align} 
\left(E_\infty\big(w,\mB_\e(x)\big)\right)^2 \ & \leq  \ \|Xu(x)\|^2 \ + \max_{\{|z-x|\leq \e\}}\Big\{
\Big[D_A \big(X_iu_\al X_iu_\al  \big)(x) \nonumber\\
& \hspace{60pt}\ - 2\de \xi_\al \big(X_iu_\al X_{iA}\big)(x)\Big]  (z-x)_A \Big\} \  + \ o(\e) ,
\end{align}
and hence
\begin{align} 
\label{2.10}
\left(E_\infty\big(w,\mB_\e(x)\big)\right)^2 \ & \leq \ \|Xu(x)\|^2\ +\ \e\Big|
D_A \big(X_iu_\al X_iu_\al  \big)(x) \nonumber\\
&\ \ \ \ - 2\de \xi_\al X_iu_\al (x)X_{iA} (x)\Big| \ + \ o(\e) , 
\end{align}
as $\e \ri 0$. Since $u$ is an Horizontal Rank-One Absolute Minimal on $\Om$,  inequalities \eqref{2.9} and  \eqref{2.10}  imply
\begin{align} \label{2.11}
0\ & \leq \ \left(E_\infty\big(w,\mB_\e(x)\big)\right)^2 \, - \,  \left(E_\infty\big(u,\mB_\e(x)\big)\right)^2 \nonumber\\
&\leq \  \e\Big(\big|D\big(X_j u_\be X_j u_\be  \big)(x)  - 2\de \xi_\al X_iu_\al (x) X_{iA} (x)\big|\\
&\ \ \ \ \ \ \ -\ \big|D\big(X_j u_\be  X_j u_\be  \big)(x)\big|\Big) \  + \ o(\e),\nonumber
\end{align}
as $\e \ri 0$. If $D\big(X_j u_\be  X_j u_\be  \big)(x) =0$, we obtain 
\begin{align}
\Big(X_i u_\al \, X_j u_\be \,  X_iX_j u_\be\Big)(x) \ &= \ X_i u_\al(x) \, X_i \, \Big(\frac{1}{2}X_ju_\be  X_j u_\be\Big)(x) \nonumber\\
 &= \ X_i u_\al(x) \, X_{iA}(x) D_A\, \Big(\frac{1}{2}X_j u_\be  X_j u_\be\Big)(x)\\
&=\ 0,\nonumber
\end{align}
as desired. If $D\big(X_j u_\be  X_j u_\be  \big)(x) \neq0$, then Taylor expansion of  the function
\[
p \ \mapsto \ \big|D\big(X_j u_\be X_j u_\be  \big)(x) +\, p\big| \, - \, \big|D\big(X_j u_\be X_j u_\be  \big)(x)  \big| 
\]
at $p_0=0$ and evaluated at $p=- 2\de \xi_\al X_iu_\al (x) X_{iA} (x)$, \eqref{2.11} implies after letting $\e \ri 0$ that
\beq \label{2.12a}
0\ \leq  \ - 2\de \xi_\al X_iu_\al (x) X_{iA} (x) \left(\frac{D_A\big(X_j u_\be X_j u_\be  \big)(x)  }{\big|D\big(X_j u_\be X_j u_\be  \big)(x)  \big|} \right) \ + \ o(\de).
\eeq
By letting $\de \ri 0$ in \eqref{2.12a} we obtain 
\begin{align}
0 \ & \geq \ \xi_\al X_i u_\al(x) \, X_{iA}(x) D_A\, \big(X_j u_\be  X_j u_\be\big)(x) \nonumber\\
 &= \ \xi_\al X_i u_\al(x) \, X_i \, \big(X_ju_\be  X_j u_\be\big)(x)\\
&=\ \xi_\al \Big(X_i u_\al \, X_j u_\be \,  X_iX_j u_\be\Big)(x) .\nonumber
\end{align}
Since $\xi$ is arbitrary, we obtain that $\big(X_i u_\al \, X_j u_\be \,  X_iX_j u_\be\big)(x)= 0$ for any $x\in \Om$. The lemma follows.               \qed

\ms

\begin{lemma} \label{l1a} If $u$ is in $\Gamma^2_{im}(\Om)^N$, solves $X_i u_\al \, X_j u_\be \, X_i X_j u_\be =0$ and $\Om$ is a Carnot-Carath\'eodory space, then $u$ is an Horizontal Rank-One Absolute Minimal on $\Om$.
\end{lemma}

Lemma \ref{l1a} is the converse of Lemma \ref{l1} and here we need the extra connectivity assumption between pairs of points in $\Om$ with curves which have horizontal velocities. This is the \emph{only} point that this assumption is needed.

\BPL \ref{l1a}. Since $X_i u_\al \, X_j u_\be \, X_i X_j u_\be =0$ on $\Om$, we have 
\beq \label{2.12b}
X_i u_\al \, X_i \Big(\frac{1}{2}X_j u_\be \, X_j u_\be \Big)\ = \ 0.
\eeq
Since $\rk(Xu)=m\leq N$, for each $x\in \Om$, the linear map $Xu(x) : H(x)\sub \R^n \larrow \R^N$ is injective and as such there exists a left inverse $(Xu(x))^{-1}$. As a result, since $\|Xu\|^2=X_j u_\be \, X_j u_\be$, we obtain
\beq \label{2.12c}
\big((Xu)^{-1}\big)_{k \al } X_i u_\al \, X_i \Big(\frac{1}{2}\|Xu\|^2 \Big)\ = \ 0,
\eeq
which implies $X_k\big(\|Xu\|^2\big)=0$ on $\Om$. Consequently, $\|Xu\|^2$ is constant along any linear combination of the vector fields $\{X_1,...,X_m\}$ on $\Om$. 

Fix $D\Subset  \Om$. Then, we have
\begin{align}
E_\infty (u, D)\ = \ \max_{\overline{D}}  \|Xu \| \ =\  \|Xu\|(\bar{x}) ,
\end{align}
for some point $\bar{x} \in \overline{D} \sub \Om$. Fix now a $g\in \Gamma^1_0(D)$ and $\xi \in \R^N$, $|\xi|=1$. We may assume $D$ is connected. Then, since $g|_{\p D}\equiv 0$, there exists an interior horizontal critical point $\bar{y}\in D$ of $g$. By using that $Xg(\bar{y})=0$, we have
\begin{align}
E_\infty (u+g\xi, D)\ 
&= \ \max_{\overline{D}}  \big\|Xu +\xi \ot Xg\big\|\\
  &\geq \  \big\|Xu(\bar{y}) +\xi \ot Xg(\bar{y})\big\| \nonumber \\
                     &= \  \big\|Xu(\bar{y}) \big\|. \nonumber
\end{align}
Since $\Om$ is a Carnot-Carath\'eodory metric space, the points $\bar{x}$, $\bar{y}$ can be connected with an admissible horizontal curve $r : [0,T] \larrow \Om $ for which
\begin{align} 
r'(t)\ & =\ a_i(t) X_i(r(t)),  \label{3.25}
\end{align}
and $ a_i(t)a_i(t) \leq  1$ for $0\leq t\leq T$, such that $r (0)= \bar{x} $ and $r (T)= \bar{y} $. Hence, by recalling that $X_i(x)=X_{iA}(x)D_A\equiv X_{iA}(x) e_A$, we have
\begin{align}
\big(E_\infty (u+g\xi, D)\big)^2 - \big(E_\infty (u, D)\big)^2\ & \geq \  \big\|Xu (\bar{y} )\big\|^2 - \big\|Xu (\bar{x} )\big\|^2 \\
          & =\  \big\|Xu (r(T))\big\|^2 - \big\|Xu (r(0))\big\|^2 \nonumber\\
&=\ \int_0^T \frac{d}{dt} \big\|Xu(r(t)) \big\|^2dt .\nonumber
\end{align}
Consequently, since $X_i\big(\|Xu\|^2\big)=0$ on $\Om$, we conclude
\begin{align}
\big(E_\infty (u+g\xi, D)\big)^2 - \big(E_\infty (u, D)\big)^2\ &  \geq \ \int_0^T \frac{d}{dt} \Big( ( X_ju_\be )(r(t))  (X_j u_\be )(r(t)) \Big)dt \nonumber\\
&=\ \int_0^T D_A\big(\|Xu\|^2 \big)(r(t)) \, r'_A(t) dt  \nonumber\\
&=\ \int_0^T D_A\big( \|Xu\|^2  \big)(r(t)) \, a_i(t) X_{iA}(r(t)) dt\\
&=\ \int_0^T X_i\big( \|Xu\|^2  \big)(r(t)) \, a_i(t)  dt \nonumber\\
&=\ 0.  \nonumber
\end{align}
Hence, $u$ is an horizontal rank-one absolute minimal and the lemma follows.
\qed

\ms

In the next two lemmas, the essential assumption needed is that the rank of $Xu$ is constant, but it is not needed the fact that it is equal to $\min\{m,N\}$.  If however $\rk(Xu)=N\leq m$, then $[Xu]^\bot\equiv 0$ and the system $[Xu]^\bot_{\al \be}X_iX_i u_\be =0$ trivializes.

\bl \label{l2} Let $u : \Om \sub \R^n \larrow \R^N$ is such that $u(\Om)$ has Horizontally $\infty$-Minimal area. Then, $u$ solves $\|Xu\|^2[Xu]_{\al \be}^\bot X_iX_i u_\be=0$ on $\Om$.
\el

\BPL \ref{l2}. Fix $x\in \Om$, $0<\e<\dist(x,\p \Om)$ and $0<\de<1$. Choose $D:= \mB_\e(x) \Subset \Om$. Fix also a vertical vector field $\nu \in \Gamma\big([Xu]^\bot,D \big)$ and an $h \in \Gamma^1\big(\bar{D} \big)$, to be specified later. We may choose $\nu$ to be a unit vector field. By differentiating along $X_i$ the equation $\nu_\al \nu_\al =1$ we obtain
\beq \label{2.12}
\nu_\al  X_i\nu_\al \ =\ 0.
\eeq
Moreover, by differentiating the identity $\nu_\al X_iu_\al= 0$ along $X_i$ (see \eqref{eq3.1}), we obtain
\beq \label{2.13}
X_i u_\al \, X_i \nu_\al \ =\  -\nu_\al \, X_iX_i u_\al .
\eeq
We set $w:=u+\de h\nu$. By using that $\nu_\al X_iu_\al =\nu_\al X_i \nu_\al = 0$, we calculate:
\begin{align}
\|Xw\|^2\ &=\  X_i(u_\al +\de h\nu_\al)\, X_i(u_\al +\de h\nu_\al)  \nonumber\\
 &=\ \Big(X_i u_\al +\de \big[\nu_\al  X_i h + hX_i\nu_\al\big]\Big) \, \Big(X_i u_\al +\de \big[\nu_\al  X_i h+ hX_i\nu_\al\big]\Big) \\
&=\  \Big(X_i u_\al +\de \big[ h X_i  \nu_\al \big]\Big) \, \Big(X_i u_\al +\de \big[h X_i  \nu_\al \big]\Big)   + \de^2 (\nu_\al \nu_\al )(X_i h X_i h) .   \nonumber
\end{align}
Hence, by using \eqref{2.13} we obtain
\begin{align} \label{2.15}
\|Xw\|^2\ &=\ \Big(X_i u_\al +\de h X_i  \nu_\al \Big) \, \Big(X_i u_\al +\de h X_i  \nu_\al \Big)   + \de^2 \|X h\|^2 \nonumber\\
            &=\ X_i u_\al  X_i u_\al  \ +\ 2\de h \big(X_i u_\al  X_i \nu_\al \big)   + \de^2 \Big( h^2 X_i\nu_\al X_i\nu_\al  + \|X h\|^2 \Big) \\
       &=\ \|Xu\|^2 \ +\ 2\de h\big(X_i u_\al  X_i \nu_\al \big) \  +\ O(\de^2)\nonumber\\
 &=\ \|Xu\|^2 \ -\ 2\de h \big(\nu_\al  X_iX_i u_\al \big) \  +\  O(\de^2). \nonumber
\end{align}
By  \eqref{2.3} and  \eqref{2.15}, we have
\begin{align} \label{2.16}
\big(E_\infty(u,\mB_\e(x))\big)^2 \ &\leq \ \big(E_\infty(u+h\nu,\mB_\e(x))\big)^2 \nonumber\\ 
&=\ \sup_{\mB_\e(x)} \|Xw\|^2  \nonumber\\
             &\leq \ \sup_{\mB_\e(x)} \|Xu\|^2\ -\ 2\de\min_{\overline{\mB_\e(x)}}\big\{h \big(\nu_\al  X_iX_i u_\al \big) \big\}\ +\ O(\de^2)\\
&= \ \big(E_\infty(u,\mB_\e(x))\big)^2\ -\ 2\de\min_{\overline{\mB_\e(x)}}\big\{h \big(\nu_\al  X_iX_i u_\al \big) \big\}\ +\ O(\de^2).  \nonumber
           \end{align}
Hence, as $\de \ri 0$, \eqref{2.16} gives
\beq \label{2.17}
\min_{\overline{\mB_\e(x)}}\big\{h \big(\nu_\al  X_iX_i u_\al \big) \big\} \ \leq \ 0.
\eeq
We now choose as $h$ the constant function $h:= \sgn \big(\big(\nu_\al  X_iX_i u_\al \big)(x)\big)$ and by \eqref{2.17} as $\e \ri 0$ we get $|(\nu_\al X_iX_i u_\al )(x)|=0$. Since $\nu$ is an arbitrary unit vertical vector field and $x$ is an arbitrary point, we get $[Xu]_{\al \be}^\bot X_iX_i u_\be=0$ on $\Om$ and the lemma follows.                              \qed

\bl \label{l3} Let $u : \Om \sub \R^n \larrow \R^N$ solve $\|Xu\|^2[Xu]_{\al \be}^\bot X_iX_i u_\be=0$ on $\Om$. Then $u(\Om)$ has Horizontally $\infty$-Minimal area.
\el

\BPL \ref{l3}. We begin with two differential identities. For any $D\Subset  \Om$, any unit vertical vector field $\nu \in \Gamma \big([Xu]^\bot,D\big)$, any $h\in \Gamma^1(\bar{D})$, $t \in \R$ and $p\geq 2$ we have 
\beq \label{2.22a}
\frac{d}{d t}\int_D \big\|X(u+t h\nu)\big\|^{p}\ = \ p\int_D\big\|X(u+t h\nu)\big\|^{p-2}X_i(u_\al +t h\nu_\al)\, X_i(h\nu_\al),
\eeq
\begin{align} \label{2.22b}
\frac{d^2}{d t^2} \int_D \big\|X(u+t h&\nu)\big\|^{p}\ = \ p\int_D 
\big\|X(u+t h\nu)\big\|^{p-2} \big\| X(h\nu)\big\|^2  \nonumber\\
&+\ p(p-2)\int_D \big\|X(u+t h\nu)\big\|^{p-4}\Big(X_i(u_\al +t h\nu_\al) \, X_i(h\nu_\al)\Big)^2.  
\end{align}
Identities \eqref{2.22a} and \eqref{2.22b}  follow by a direct calculation and by using that $\|Xu\|^2=X_iu_\al X_iu_\al $. We set:
\[
f(t)\ := \ \int_D \big\|X(u+t h\nu)\big\|^{p} - \int_D \big\|Xu\big\|^{p}
\]
Evidently, $f$  vanishes at $t=0$. By \eqref{2.22b}, $f$ is convex. Moreover, we have
\begin{align} \label{2.23a}
f'(0)\ &=\ \frac{d}{d t}\Big|_{t=0} \int_D \big\|X(u+th\nu)\big\|^{p} \\
 &=\ p\int_D \|Xu\|^{p-2}X_i u_\al X_i (h\nu_\al), \nonumber
\end{align} 
and since $\nu_\al X_i u_\al=0$, we obtain
\begin{align} \label{2.23a}
f'(0)\ &= \ p\int_D \|Xu\|^{p-2} X_i u_\al \big( h X_i\nu_\al\, +\, \nu_\al X_ih\big)  \nonumber\\
&= \ p\int_D \|Xu\|^{p-2} X_i u_\al \big( h X_i\nu_\al\big)\\
&= \ -p\int_D  \|Xu\|^{p-2} h(\nu_\al X_iX_i u_\al ).  \nonumber
\end{align} 
Since $u$ solves $\|Xu\|^2[Xu]_{\al \be}^\bot X_iX_iu_\be=0$ with $\|Xu\|>0$ on $\Om$, we deduce that $f'(0)=0$. Since $f$ is convex, $t=0$ is a point of global minimum for $f$ and hence $f(t)\geq f(0)$. Thus,
\beq \label{2.24a}
\int_D \big\|Xu\big\|^{p} \ \leq \int_D \big\|X(u+t h\nu)\big\|^{p},
\eeq
for any $t\in \R $. By rescaling \eqref{2.24a} and letting $p\ri \infty$ we conclude that 
\[
E_\infty(u,D)\ \leq\ E_\infty(u+h\nu,D)
\]
and consequently $u(\Om)$ has Horizontaly $\infty$-Minimal area on $\Om$. The lemma has been established.            \qed

\ms
In view of Lemmas \ref{l1}, \ref{l1a}, \ref{l2} and \ref{l3}, Theorem \ref{th1} follows.

\section{Maximum and Minimum Principles for $\|Xu\|$ for the subelliptic $\infty$-Laplace System.} 

In this brief section we establish maximum and minimum principles for the induced (from the Riemannian metric) norm of the horizontal gradient of solutions of full horizontal rank to the subellptic $\infty$-Laplacian on Carnot-Carath\'eodory spaces. Herein we subsume the same setting introduced and utilized in the previous sections.

\begin{proposition}[Maximum-Minimum Principles for the Horizontal Gradient on Carnot-Carath\'eodory spaces] Suppose \label{pr1} $u :\Om \sub \R^n \larrow \R^N$ is in $\Gamma^2_{rk}(\Om)^N$, $\De^X_\infty u =0$ and assume $\Om$ is a Carnot-Carath\'eodory space. Then, for any $D \Subset \Om$ we have:
 \begin{align} 
\sup_{D}\|Xu\|\ = \ \max_{\p D}\|Xu\|, \label{2.23}\\
\inf_{D}\|Xu\|\ = \ \min_{\p D}\|Xu\|. \label{2.24}
\end{align}
\end{proposition}

The connectivity assumption by admissible curves is needed in the case of horizontal immersions. The main tool in the case of horizontal submersions is the following

\begin{lemma}[Subelliptic gradient flow] \label{l5}
Let $u :\Om \sub \R^n \larrow \R^N$ be in $\Gamma^2(\Om)^N$. Consider the gradient flow
\beq \label{2.26}
\left\{
\begin{array}{l}
r'_A(t)\ = \ \left(      \dfrac{\|Xu\|^2}{\xi_\al X_iu_\al \xi_\be X_i u_\be }  
X_{iA} \xi_\ga X_i u_\ga\right)\big(r(t)\big), \ \ t\neq 0,\ms\\
r(0)\ = \ x,
\end{array}
\right.
\eeq
for $x\in \Om$, $\xi \in \R^N$, $|\xi|=1$ and $\xi_\al Xu_\al(x) \neq 0$, where $X_i(x)=X_{iA}(x)D_A$ in the standard frame. Then, we have the differential identities
\begin{align}
\ \ \ \ \ \ \ \  \frac{d}{d t} \Big(\dfrac{1}{2}\big\|Xu\big(r(t)\big)\big\|^2\Big)\ &= \ \left(  \dfrac{\|Xu\|^2}{\xi_\al X_iu_\al \xi_\be X_iu_\be } \xi_\ga \big[X_i u_\ga X_j u_\de X_iX_ju_\de \big] \right)\big(r(t)\big), \label{2.27}\\
\frac{d}{d t}\Big(\xi_\al u_\al \big(r(t)\big)\Big)\ &= \ \big\| Xu \big(r(t)\big)\big\|^2, \label{2.28}
\end{align}
which imply $X_i u_\al X_j u_\be X_iX_ju_\be=0$ on $\Om$ if and only if for each $x\in 
\Om$ and every parameter $\xi$, $\|Xu\big(r(t)\big)\|$ is constant along the trajectory $r(t)$ and $t\mapsto \xi_\al u_\al \big(r(t)\big)$ is \emph{affine}.
\end{lemma}

For the shake of completeness, let us give the elementary proof of this lemma.

\BPL \ref{l5}. Using the flow \eqref{2.26}, we calculate
\begin{align}
\frac{d}{d t} \Big(\dfrac{1}{2}\big\|Xu\big(r(t)\big)\big\|^2\Big)\ &=\ \Big(X_ju_\de D_A X_ju_\de \Big)\Big|_{r(t)}r'_A(t) \nonumber\\
 &=\ \Big(X_ju_\de D_A X_ju_\de \Big)\Big|_{r(t)}\left(      \dfrac{\|Xu\|^2}{\xi_\al X_iu_\al \xi_\be X_i u_\be }   X_{iA} \xi_\ga X_i u_\ga\right) \Big|_{r(t)}\\
&= \ \left(  \dfrac{\|Xu\|^2}{\xi_\al X_iu_\al \xi_\be X_iu_\be } \xi_\ga \big[X_i u_\ga X_j u_\de X_iX_ju_\de \big] \right)\Big|_{r(t)},  \nonumber
\end{align}
and also
\begin{align}
\frac{d}{d t}\Big(\xi_\ga u_\ga \big(r(t)\big)\Big)\ &=\ \xi_\ga D_Au_\ga(r(t))r'_A(t)  \nonumber\\
&=\ \left(      \dfrac{\|Xu\|^2}{\xi_\al X_iu_\al \xi_\be X_i u_\be }  
 \xi_\ga D_Au_\ga \, \xi_\de \, X_{iA}  X_i u_\de  \right) \Big|_{r(t)}\\
&=\ \left(      \dfrac{\|Xu\|^2}{\xi_\al X_iu_\al \xi_\be X_i u_\be }  
 \xi_\ga X_iu_\ga \, \xi_\de  X_i u_\de  \right) \Big|_{r(t)}   \nonumber\\
&= \ \big\| Xu \big(r(t)\big)\big\|^2  \nonumber
\end{align}
The lemma readily follows. \qed

\BPP \ref{pr1}.  Fix $D \Subset \Om$. Consider first the case of horizontal immersions where $\rk(Xu)=m\leq N$. We argue as in Lemma \ref{l1a}: by \eqref{eq2} we have $X_i u_\al X_i \big(\frac{1}{2}\|Xu\|^2\big) =0$. Since for each $y \in \Om$ the linear map $Xu(y) : H(y) \sub \R^n \larrow \R^N$ is injective, $Xu$ is left-invertible and we obtain $X \big(\frac{1}{2}\|Xu\|^2\big) =0$ on $\Om$. We also have 
\[
\sup_D\|Xu\|\ =\ \|Xu(\bar{x})\|\ , \ \ \ \max_{\p D} \|Xu\|\ =\ \|Xu(\bar{y})\|, 
\]
for some points $\bar{x}, \bar{y} \in \Om$. Since $\Om$ is a Carnot-Carath\'eodory space, we can connect them with an admissible curve $r: [0,T]\larrow \Om$ for which $r(0)=\bar{x}$ and $r(T)=\bar{y}$ with horizontal velocity: $r'(t) \in H(r(t))$. Then,
\begin{align}
\max_{\p D} \|Xu\|^2 - \sup_D\|Xu\|^2 \ &=\   \|Xu(\bar{y})\|^2 -  \|Xu(\bar{x})\|^2\nonumber\\
&=\ \int_0^T \frac{d}{dt} \|Xu(r(t))\|^2dt\\
&=\ 2 \int_0^T \big( X_j u_\be  D_A X_j u_\be \big) (r(t)) r'_A(t)dt\nonumber,
\end{align}
and hence, for some coefficients $a_i(t)$ (since $X_i=X_{iA}D_A$),
\begin{align}
\max_{\p D} \|Xu\|^2 - \sup_D\|Xu\|^2 \
&=\ 2 \int_0^T \big( X_j u_\be  D_A X_j u_\be \big) (r(t))\, a_i(t)X_{iA}(r(t)) dt\nonumber\\
&=\ 2 \int_0^T \big(X_j u_\be  X_i X_j u_\be \big)(r(t))\, a_i(t)dt\\
&=\  \int_0^T X_i \big(X_j u_\be X_j u_\be \big)(r(t)) \, a_i(t)dt\nonumber\\
&=\ 0.\nonumber
\end{align}
Hence, \eqref{2.23} follows and the same holds for \eqref{2.24}.

For the case of horizontal submersions where $\rk(Xu)=N\leq m$, fix $x\in D\Subset \Om$ such that
\beq
\sup_{D}\|Xu\| \ = \ \|Xu(x)\|  .
\eeq
Choose also $\xi \in \R^N$ with $|\xi|=1$ and consider the gradient flow \eqref{2.26} at $x$. Since $\rk(Xu)=N\leq m$, for each $y\in \Om$ the linear mapping $Xu(y) : H(y)\sub \R^n \larrow \R^N$ is surjective and hence 
\[
\big(\xi_\al X_i u_\al \, \xi_\be X_i u_\be\big)(y)\ >\ 0
\]
for all $y \in \Om$. Hence, the flow is globally defined on $\Om$ for all parameters $\xi$.  By \eqref{2.27}, we have
\[
\|Xu\big(r(t)\big)\|\ =\ \|Xu(x)\| 
\]
along the trajectory and by utilizing \eqref{2.28} we see that the trajectory reaches $\p D$ in finite time (since $D$ is bounded) because
\begin{align}
\xi_\al u_\al ( r(t))\ - \ \xi_\al u_\al (x)\ =\  t\|Xu(x)\|^2.
\end{align}
Hence, there exists $T=T(x)>0$ such that $r(T) \in \p D$. Consequently,
\begin{align}
\sup_{D}\|Xu\| \ &= \ \|Xu(x)\|  \nonumber\\
   &= \ \|Xu\big(r(T)\big)\|\\
& \leq \  \max_{\p D} \|Xu\|\nonumber
\end{align}
and similarly we obtain $\inf_D\|Xu\|\geq \min_{\p D}\|Xu\|$. The proposition readily follows, in view of the observation that $\sup_D \|Xu\|=\max_{\overline{D}} \|Xu\|$ because $\|Xu\| \in C^0(\overline{D})$ and $D \Subset \Om$.                      \qed

\ms

\ms

\noi \textbf{Acknowledgement.} The author has greatly benefited from the scientific discussions with Sergio Polidoro during his visit at BCAM, whom he would like to deeply thank for selflessly sharing his expertise on subelliptic equations. He is also indebted to the anonymous referee for the careful reading of the manuscript and for his suggestions which improved the content and the presentation of the paper.

\ms

\bibliographystyle{amsplain}

\end{document}